\documentclass[12pt,a4paper]{article}
\usepackage[T1]{fontenc}
\usepackage[utf8]{inputenc}
\usepackage{graphicx}
\usepackage{amscd}
\usepackage{amsmath,amssymb}

\usepackage{theorem,url}
\usepackage{subfig}
\usepackage[version=3]{mhchem}

\renewcommand{\emph}{\textbf}
\newcommand{\mma}[1]{\textbf{\texttt{#1}}}

\theorembodyfont{\upshape}
\newtheorem{Thm}{Theorem}
\newtheorem{Def}{Definition}
\newtheorem{Exa}{Example}
\newtheorem{Rem}{Remark}

\newtheorem{Cor}{Corollary}
\newtheorem{Sta}{Statement}

\newcommand{\torol}[1]{}

\newcommand{\nulb}{\ifmmode \mathbf{0}\else \textbf{0}\fi}

\newcommand{\cb}{\mathbf{c}}

\newcommand{\xb}{\mathbf{x}}
\newcommand{\N}{\mathbb{N}}
\newcommand{\Pb}{\mathbf{P}}
\newcommand{\R}{\mathbb{R}}
\newcommand{\Rp}{\mathbb{R}^+}
\newcommand{\SSS}{\mathcal{S}}

\newcommand{\ikde}{induced kinetic differential equation}

\newcommand{\spanlin}{\mathop{\mathrm{span}}\nolimits}
\newcommand{\Mma}{\textit{Mathematica }}

\title{Quadratic First Integrals of\\ Kinetic Differential Equations}
\author{
Ilona Nagy\\
Department of Mathematical Analysis\\
Budapest University of Technology and Economics,\\
Budapest, Egry J. u. 1., HUNGARY, H-1111\\
\url{nagyi@math.bme.hu} and\\
J\'anos T\'oth\\
Department of Mathematical Analysis\\
Budapest University of Technology and Economics,\\
Budapest, Egry J. u. 1., HUNGARY, H-1111 and\\
Laboratory for Chemical Kinetics
Eötvös University\\
Budapest, Pázmány Péter sétány 1/A, HUNGARY, H-1117\\
\url{jtoth@math.bme.hu}}
\date{}

\begin{document}
{\LARGE \noindent \textbf{Title Page}

\noindent \textbf{Title:} Quadratic First Integrals of Kinetic
Differential Equations

\noindent \textbf{Authors:}

Ilona Nagy\footnote{(author designated to review proofs) Department
of Mathematical Analysis, Budapest University of Technology and
Economics, Budapest, Egry J. u. 1., HUNGARY, H-1111, Phone: 361
463-5141, Fax: 361 463 3172, E-mail: \url{nagyi@math.bme.hu}}

Tóth, J.\footnote{ Department of Analysis, Budapest University of
Technology and Economics, Egry J. u. 1., Budapest, Hungary, H-1111
and Laboratory for Chemical Kinetics, Eötvös Loránd University,
Pázmány P. sétány 1/A, Budapest, Hungary, H-1117}

\noindent \textbf{Running Head}: First integrals in kinetics}
\vfill\eject

\maketitle

\begin{abstract}
Classes of kinetic differential equations are delineated which do
have a quadratic first integral, and classes which can not have one.
Example reactions corresponding to the obtained kinetic differential
equations are shown, and a few figures showing the trajectories of
the corresponding systems are also included. Connections to other
areas are mentioned and unsolved problems collected. The new results
are theoretical, although computational tools are heavily used.
Applications from biology and combustion theory will come later.
\end{abstract}

\section{Introduction}

Our aim is to determine classes of mass action type kinetic
differential equations with the property of having a quadratic first
integral. Since the introduction of the name of first integral by E.
N\"other, in 1918 it turned out that first integrals may help
\begin{itemize}
\item
prove that the complete solution of the induced kinetic differential
equation is defined for all positive times \cite[p. 586, Theorem
9]{volpert};
\item
reduce the number of variables either by constructing an appropriate
lumping scheme \cite{lirabitztoth} or by simply eliminating some
variables;
\item
apply the generalization of the Bendixson and Bendixson--Dulac
criterion to higher dimensional cases
\cite{tothZAMM,webersturmseilerabdelrahman}.
\end{itemize}

Our main tool to find such first integrals is the comparison of
coefficients of polynomials and the characterization of kinetic
differential equations within the class of polynomial ones
\cite{harstoth}. (A good review of our earlier results has been
given by \cite{chellaboinabhathaddadbernstein}.) This
characterization has proved quite useful in
\begin{itemize}
\item
designing minimal oscillatory reactions \cite{tothharsTCA},
\item
providing an alternative proof for the uniqueness of the
Lotka--Volterra model \cite{pota,schumantoth},
\item
investigating chaos in chemical reactions
\cite{tothharsPhys,halmschlagerszenthetoth},
\item
investigating symmetries in kinetic differential equations
\cite{tothgradient,totherdisymmetries};
\item
selecting the kinetic lumping schemes from all the possible ones
\cite{farkaslumping},
\item
finding necessary conditions of Turing instability
\cite{szilitothPRE,szilitothJMC}.
\end{itemize}

(Let us remark in passing that Dilao \cite{dilao} in his detailed
analysis of Turing instability disregarded this characterization,
therefore his Case f) in Theorem 2.2 cannot occur in a mass action
type kinetic model.)

We are also interested in kinetic differential equations with
quadratic first integrals which describe mass conserving reactions
\cite{hornjackson,deaktothvizvari,schusterhofer}. It turns out in
some cases that the existence of a quadratic first integral and mass
conservation together form a too rigorous set of requirements: we
may be able to prove that such equations do not exist.

Once we have a kinetic differential equation fulfilling some
requirements we might be interested in reactions with the given
induced kinetic differential equation. However, the solution to this
problem is far from being unique \cite[page 48--49]{tothsztaki},
\cite[page 67--68]{erditoth}. One possible approach might be that we
try to find a reaction with a given property (weak reversibility,
zero deficiency etc.), or a minimal or maximal reaction in a certain
sense with a given property \cite{szederkenyi,szederkenyihangos}.

The structure of our paper is as follows. Section 2 presents the
basic definitions. Section 3 gives the general results, both
positive and negative: on the existence and nonexistence of kinetic
differential equations with quadratic first integrals depending
possibly on further assumptions. In some cases we also show a
reaction having the obtained differential equation, and a few
figures reflecting the behavior of possible trajectories. Section 4
shows how \Mma (more precisely, the program package
\mma{ReactionKinetics} written in the \Mma\ language,
\cite{nagypapptoth,tothnagypapp, tothnagyzsely}) can be used to
formulate, prove or disprove conjectures. Finally, Section 5
formulates problems to be solved.
\section{MASS ACTION TYPE KINETIC DIFFERENTIAL EQUATIONS}
Here we recapitulate very shortly the basic concepts of formal
reaction kinetics as they can be found e.g. in \cite{erditoth},
\cite{yablonskymarin} or \cite{tothnagypapp}.
\subsection{Induced kinetic differential equation of a reaction with mass action type kinetics}
Let us consider a vessel (a cell, a reactor, a test tube etc.) of
constant volume at constant pressure and temperature and let
\(M,R\in\N; \alpha,\beta\in\N_0^{M\times R},\) and consider the
\emph{complex chemical reaction}:
\begin{equation}\label{ccr}
\ce{\sum_{m=1}^M \alpha(m,r)X(m) -> \sum_{m=1}^M
\beta(m,r)X(m)}\quad (r=1,2,\dots,R),
\end{equation}
where the components of the matrices
\(\alpha=(\alpha(m,r))_{m=1,2,\dots,M;r=1,2,\dots,R}\) and
\(\beta=(\beta(m,r))_{m=1,2,\dots,M;r=1,2,\dots,R}\) are the
\emph{stochiometric coefficients}. There are a few natural
conditions fulfilled by the stoichiometric matrix (see e.g.
\cite[page 77]{deaktothvizvari}):
\begin{enumerate}
\item all the species take part in at least one reaction step;
\item all the reaction steps change the quantity of at least one species;
\item all the reaction steps are determined by their
reactant and complex products.
\end{enumerate}
\begin{Rem}
The last requirement may be too restrictive because if a reaction
can proceed through two different transition states, as in the
reaction
\begin{equation*}
\ce{CH3CHOH + O2 <=> HO2 + CH3CHO}
\end{equation*}
\cite{zadorfernandesgeorgievskiimelonitaatjesmiller}, then one
should duplicate this step to exactly represent the mechanistic
details of the reaction.
\end{Rem}

Suppose the reaction can adequately be described using \emph{mass
action kinetics}, then its \emph{deterministic model} is
\begin{eqnarray}\label{ikde}
{c_m}'(t)&=&f_m(\cb(t)):=
\sum_{r=1}^R(\beta(m,r)-\alpha(m,r))k_r\prod_{p=1}^{M}c_p(t)^{\alpha(p,r)}\\
c_m(0)&=&c_{m0}\in\R_0^+\quad(m=1,2,\dots,M)\label{ini}
\end{eqnarray}
(with the positive \emph{reaction rate coefficients}
$k_r$)---describing the time evolution of the concentration vs. time
functions
$$
t\mapsto c_m(t):=[X(m)](t)
$$
of the species. Eq. \eqref{ikde} is also called the mass action type
\emph{induced kinetic differential equation} of the reaction
\eqref{ccr}.

\subsection{Polynomial and kinetic differential equations}
The induced kinetic differential equation of the reaction
\eqref{ccr} is a polynomial differential equation, because all the
functions $f_m$ are polynomials in all their variables. (This
property can be shown to be equivalent with the fact that $f_m$ is
an $M$-variable polynomial \cite{carroll}.) However, it is not true
that all polynomial differential equations can be obtained as
induced kinetic differential equations of some reactions, as the
examples
$$
x'=y,\quad y'=\boxed{-x}
$$
(of the harmonic oscillator), or the Lorenz model
$$
x'=\sigma(y-x),\quad y'=\rho x-\boxed{xz}\, ,\quad z'=xy-\beta
z\quad (\sigma,\rho,\beta>0)
$$
show. The speciality of kinetic differential equations is that they
cannot contain terms like those boxed above, i.e. terms expressing
the decay of a quantity without its participation. Such terms are
said to represent \emph{negative cross effects} \cite{harstoth}.
Moreover, it is also true that the absence of such terms allows us
to construct a reaction inducing the given differential equation
\cite{harstoth}. To formulate this property and also our statements
below we need the following definition.
\begin{Def}
Let \(M\in\N,\) and let us suppose that
\(\Pb:\R^M\longrightarrow\R^M\) is a function with the property that
all its coordinate functions are polynomials in all their variables.
Then the differential equation
\begin{equation}
{\xb}'=\Pb\circ\xb\label{polydeq}
\end{equation}
is said to be a \emph{polynomial differential equation}.
\end{Def}

Let us remark that $\mathbb{R}^2\ni(x,y)\rightarrow xy$ is a second
degree polynomial although it is of the first degree in all of its
variables. When formulating and proving our results in Section 3
below we sometimes need notations different from those in the above
definition for the sake of transparency.

\begin{Def}
Let us consider the polynomial differential equation
\eqref{polydeq}, and suppose that there is an
\(m\in\{1,2,\dots,M\}\) and a vector
$$\cb^m_0:=(c_1,c_2,\dots,c_{m-1},0,c_{m+1},\dots,c_M)$$ $(c_p\ge0 \text{\ for\ }
p=1,2,\dots,m-1,m+1,\dots,M)$ so that \(P_m(\cb^m_0)<0\). Then,
\eqref{polydeq} is said to contain \emph{negative cross-effect}.
\end{Def}

Our starting point is the following statement the constructive proof
of which can be found in \cite{harstoth}.
\begin{Thm}
A polynomial differential equation is the induced kinetic
differential equation of a reaction endowed with mass action type
kinetics if and only if it contains no negative cross-effect.
\end{Thm}
Then, it is quite natural to call polynomial differential equations
\emph{kinetic} if they have no negative cross-effect.

Let us remark that the absence of a negative cross-effect is
stronger than the property that the velocity of the vector field is
always pointing into the interior of the first orthant (implying
that the first orthant is an invariant set of \eqref{ikde}) as the
remark by Feinberg (cited in \cite[p. 41]{tothsztaki}) shows:
\begin{equation}
{c}_1'=c_2+c_2^2-2c_2c_3+c_3^2,\quad {c}_2'=0,\quad
{c}_3'=0.\nonumber
\end{equation}
It may be useful to know that Chellaboina et al. gave practically
the same example and also reproduced our proof of the above theorem
in \cite{chellaboinabhathaddadbernstein}.

\subsection{Mass conservation}
Although it is very convenient to allow reactions like \ce{X -> 0}
to describe outflow, or those like \ce{0 -> X} to represent inflow,
or \ce{X -> 2X} to denote autocatalysis, it is still quite natural
to give extra importance to reactions which do conserve mass. The
intuitive meaning of mass conservation is that calculating the total
mass on both sides of a reaction step we get the same amount
\cite[page 89]{hornjackson}.
\begin{Def}
The reaction \eqref{ccr} is said to be \emph{stoichiometrically mass
conserving}, if there exists a vector $\varrho\in(\Rp)^M$ for which
\begin{equation}\label{mass}
\forall
r\in\{1,2,\dots,R\}:\sum_{m=1}^M\varrho(m)\alpha(m,r)=\sum_{m=1}^M\varrho(m)\beta(m,r)
\end{equation}
is fulfilled.
\end{Def}
It is not so trivial to decide if a reaction of the form \eqref{ccr}
is stoichiometrically mass conserving or not if we are only given
the stoichiometric coefficients
\cite{deaktothvizvari,schusterhofer}. (These last papers provide
sufficient and necessary conditions of, and algorithms to decide
mass conservativity.)

Now an equivalent definition of stoichiometric mass conservation
will be given. To arrive at that definition preparations are to be
made.

\begin{Def}
The \emph{stochiometric subspace} of the reaction \eqref{ccr} is the
linear space
$$
\SSS:=\spanlin{\{\alpha(.,r)-\beta(.,r); r=1,2,\dots,R\}}.
$$
With this, the reaction \eqref{ccr} is stoichiometrically mass
conserving if there exists a vector with positive coordinates in the
orthogonal complement of the stoichiometric subspace. The set
$c_0+\SSS\quad(c_0\in(\Rp)^M)$ is a (positive) \emph{reaction
simplex}.
\end{Def}
A fundamental result by Horn and Jackson \cite{hornjackson} follows.
\begin{Thm}
A reaction is stoichiometrically mass conserving if and only if all
positive reaction simplexes are bounded.
\end{Thm}

An immediate consequence of the theorem is that a complete solution
of a stoichiometrically mass conserving system is defined for all
nonnegative times. Further much more refined statements on
nonnegativity can be found in \cite{volpert}.
\begin{Exa}
Stoichiometric mass conservation is sufficient but not necessary for
the preservation of (possibly, weighted) total mass. The example
\cite[page 89]{feinberghorn}
\begin{figure}[!ht]
\includegraphics[height=2cm]{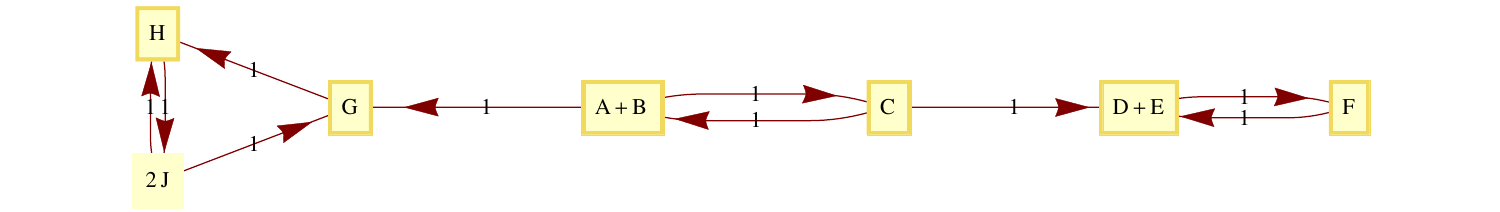}
\caption{The Feinberg--Horn example{\label{fig:1}}}
\end{figure}
shows that a positive linear combination of the concentrations can
be constant even if the positive coefficients do not lie in the
orthogonal complement of the stoichiometric subspace. If this is the
case one may speak about a \emph{kinetically mass conserving}
reaction. Let us see the details.

Let all the reaction rate constants be unity. Now we are going to
show that the vector
$\varrho=\left(\begin{array}{lllllllll}1&2&4&1&4&5&2&2&1\end{array}\right)^{\top}$
is orthogonal to the right hand side of the \ikde\
\begin{align*}
a'&=-2ab+c   &b'&=-2ab+c   &c'&=ab-2c\\
d'&=c-de+f   &e'&=c-de+f   &f'&=de-f\\
g'&=ab-g+j^2  &h'&=g-h+j^2  &j'&=2h-4j^2
\end{align*}
of the reaction but it is not orthogonal from the left to the matrix
\begin{equation}
\gamma= \left(
\begin{array}{r|rrrrrrrrrr}
 &1&2&3&4&5&6&7&8&9&10\\
      \hline
\ce{A}& -1 & 1 & 0 & 0 & 0 & -1 & 0 & 0 & 0 & 0 \\
\ce{B}&-1 & 1 & 0 & 0 & 0 & -1 & 0 & 0 & 0 & 0 \\
\ce{C}& 1 & -1 & -1 & 0 & 0 & 0 & 0 & 0 & 0 & 0 \\
\ce{D}& 0 & 0 & 1 & -1 & 1 & 0 & 0 & 0 & 0 & 0 \\
\ce{E}& 0 & 0 & 1 & -1 & 1 & 0 & 0 & 0 & 0 & 0 \\
\ce{F}& 0 & 0 & 0 & 1 & -1 & 0 & 0 & 0 & 0 & 0 \\
\ce{G}& 0 & 0 & 0 & 0 & 0 & 1 & -1 & 0 & 0 & 1 \\
\ce{H}& 0 & 0 & 0 & 0 & 0 & 0 & 1 & -1 & 1 & 0 \\
\ce{J}& 0 & 0 & 0 & 0 & 0 & 0 & 0 & 2 & -2 & -2
\end{array}
\right)
\end{equation} of the elementary reaction vectors. Really,
$$
\varrho^{\top}\left(\begin{array}{lllllllll}a'&b'&c'&d'&e'&f'&g'&h'&j'\end{array}\right)^{\top}=0,
$$
and
$$
\varrho^{\top}\gamma=\left(\begin{array}{llllllllll}1&-1&1&0&0&-1&0&0&0&0\end{array}\right)^{\top}\neq\nulb^{\top}.
$$

\end{Exa}
It is a very natural requirement that a numerical method aimed at
solving \eqref{ikde} should keep the total mass
$\sum_{m=1}^{M}\varrho_mc_m(t)$ constant (independent on time) in
case of a mass conserving reaction. There are some methods to have
this property, see e.g. \cite{bertolazzi}. A similar requirement is
to keep other, e.g. quadratic first integrals, what has also been
shown for some methods \cite{rosenbaum}.

However, not much is known about equations, especially kinetic
differential equations with quadratic first integrals. Obviously,
equations of mechanics, like that of the standard harmonic
oscillator $ x'=y\quad y'=-x $ may have quadratic first integrals,
$V(p,q):=p^2+q^2$ in this case, and here the meaning of the
quadratic first integral is the total mechanical energy.
\section{EXISTENCE AND NONEXISTENCE OF\\
 QUADRATIC FIRST INTEGRALS}
\subsection{Diagonal first integrals}

\begin{Thm}\label{th:elsotetel}
Let us consider the following system of differential equations
\begin{eqnarray}
x_m'=F_m\circ(x_1,x_2,\dots,x_M), \quad (m=1,\dots,M)\label{system1}
\end{eqnarray}
where the functions $F_m$ are quadratic functions of the variables,
that is,

\begin{eqnarray}\label{quadratic1}
F_m(x_1,x_2,\dots,x_M) =\sum_{p=1}^{M}A_{m,p}x_p^2
+\sum_{\begin{smallmatrix}p=1\\p\neq m\end{smallmatrix}}^{M}B_{m,p}x_m x_p\nonumber\\
+\sum_{\begin{smallmatrix}p,q=1\\p<q\\p\neq m,q\neq
m\end{smallmatrix}}^{M}C_{p,q}^m x_p x_q +\sum_{p=1}^M
D_{m,p}x_p+E_m.
\end{eqnarray}
Suppose that the system of differential equations is \emph{kinetic}.
The function
$$
V(x_1,x_2,\dots,x_M)=a_1x_1^2+a_2x_2^2+\dots +a_Mx_M^2$$ (with
$a_m>0$ for $m=1,2,\dots,M$) is a first integral for the above
system if and only if the functions $F_m$ have the following form
with $K_{m,p}\geq0$:
\begin{equation}\label{ered:elsotetel}
F_m(x_1,x_2,\dots,x_M) =\sum_{\begin{smallmatrix}p=1\\p\neq
m\end{smallmatrix}}^{M}a_p K_{m,p}x_p^2
-\sum_{\begin{smallmatrix}p=1\\p\neq m\end{smallmatrix}}^{M}a_p
K_{p,m}x_m x_p.
\end{equation}
\end{Thm}

\emph{Proof}. The function $V$ is a first integral for the system
\eqref{system1} if and only if its Lie-derivative with respect to
the system is equal to zero, that is,
\begin{eqnarray}
0&=&\frac{1}{2}\sum_{m=1}^M \partial_m V(x_1,x_2,\dots,x_M)F_m(x_1,x_2,\dots,x_M)\nonumber\\
&=&\sum_{m=1}^M a_m\left(\sum_{p=1}^M A_{m,p}x_m x_p^2+
\sum_{\begin{smallmatrix}p=1\\p\neq m\end{smallmatrix}}^{M}B_{m,p}x_m^2 x_p\right.\nonumber\\
&+&\left.\sum_{\begin{smallmatrix}m\neq p,m\neq
q\\p<q\end{smallmatrix}}C_{p,q}^m x_m x_p x_q+ \sum_{p=1}^M
D_{m,p}x_m x_p+E_m x_m\right)\label{Lie-derivative}
\end{eqnarray}

Since the system \eqref{system1} is kinetic, the coefficients of
terms in $F_m$ not containing $x_m$ are nonnegative:

\begin{tabular}{lllllllll}
$(k1)$& $A_{m,p}$&$\ge$&$0$&for&$m,p$&$=1,2,\dots M;$&$ m\neq p$\\
$(k2)$& $C^m_{p,q}$&$\ge$&$0$&for&$m,p,q$&$=1,2,\dots M;$&$ m\neq p,m\neq q, p<q$\\
$(k3)$& $D_{m,p}$&$\ge$&$0$&for&$m,p$&$=1,2,\dots M;$&$ m\neq p$\\
$(k4)$& $E_m$&$\ge$&$0$&for&$m$&$=1,2,\dots M.$&$$
\end{tabular}

For all $m$, the coefficients of $x_m^3$, $x_m^2$ and $x_m$ are
$a_mA_{m,m}$, $a_mD_{m,m}$ and $a_mE_{m,m}$, respectively. Since
these monomials  are independent of each other and of the other
terms in \eqref{Lie-derivative}, it follows that
$A_{m,m}=D_{m,m}=E_{m,m}=0$.

If $m\neq p$, then the monomial $x_m x_p$ appears twice in
\eqref{Lie-derivative} with coefficients $a_m D_{m,p}$ and $a_p
D_{p,m}$. Thus $a_m D_{m,p}+a_p D_{p,m}=0$ and because of $(k3)$,
$D_{m,p}=D_{p,m}=0$.

If $m\neq p,m\neq q,p<q$, then the monomial $x_m x_p x_q$ appears
three times in \eqref{Lie-derivative} with coefficients $a_m
C_{p,q}^m$, $a_p C_{m,q}^p$ and $a_q C_{m,p}^q$. Thus $a_m
C_{p,q}^m+a_p C_{m,q}^p+a_q C_{m,p}^q=0$ and because of $(k2)$,
$C_{p,q}^m=C_{m,q}^p=C_{m,p}^q=0$.

If $m\neq p$, then the monomial $x_m x_p^2$ appears twice in
\eqref{Lie-derivative} with coefficients $a_m A_{m,p}$ and $a_p
B_{p,m}$ and thus $a_m A_{m,p}+a_p B_{p,m}=0$ where $A_{m,p}\geq0$
because of $(k1)$. Without the loss of generality, it may be assumed
that $A_{m,p}=a_p K_{m,p}$ where $K_{m,p}\geq0$ and so $B_{p,m}=-a_m
K_{m,p}$.

The the proof of the \emph{if} part is obvious.

\begin{figure}[h!]
  \centering
  \subfloat[][]{\label{fig:3a}\includegraphics[width=0.4\textwidth]{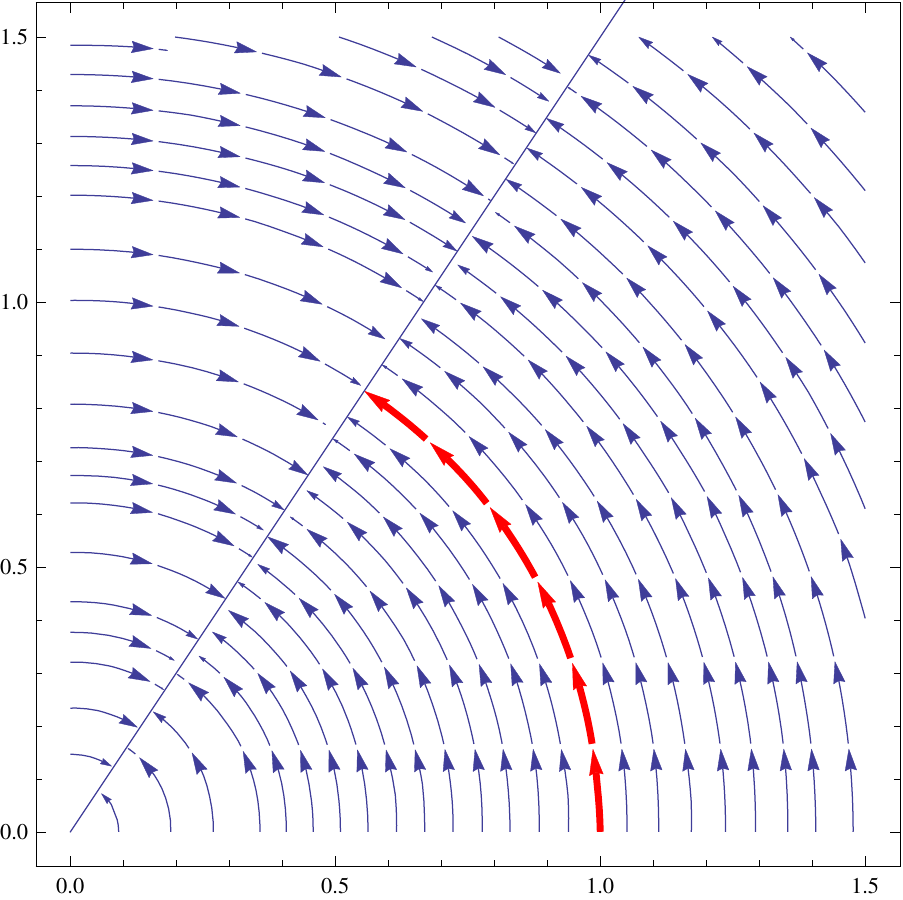}}\quad
  \subfloat[][]{\label{fig:3b}\includegraphics[width=0.4\textwidth]{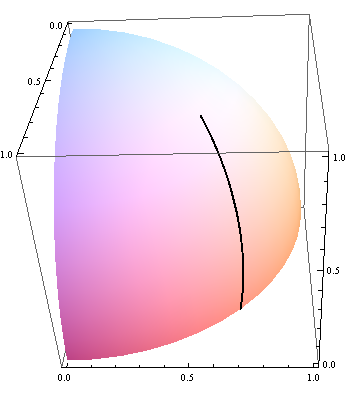}}
  \caption{(a) Trajectories of system \eqref{exa:elsotetel2D} with $a=2, b=3$ starting from
$x(0)=1, y(0)=0.$\\ (b) Trajectories of system
\eqref{exa:elsotetel3D} with $a=2, b=3, c=4,d=5,e=6,f=7$ starting
from $x(0)=\frac{1}{\sqrt{2}}, y(0)=\frac{1}{\sqrt{2}},
z(0)=0.$}\label{fig:3}
\end{figure}

\begin{Exa}
Let $M=2$ and suppose that $V(x,y)=x^2+y^2$. Then
\eqref{ered:elsotetel} specializes to
\begin{equation}\label{exa:elsotetel2D}
{x}'=ay^2-bxy,\quad{y}'=bx^2-axy
\end{equation}
which may be considered as the \ikde\ of the reaction
\begin{equation}
\ce{X <-[a] X + Y ->[b] Y}\qquad\ce{2X ->[b] 2X + Y}\qquad\ce{2Y
->[a] X + 2Y}
\end{equation}
as the application of \mma{RightHandSide[\{X <- X + Y -> Y, 2 X -> 2
X + Y,
  2 Y -> X + 2 Y\}, \{a, b, b, a\}, \{x, y\}]}
gives: $\{a y^2-b x y,b x^2-a x y\}.$ A typical trajectory is shown
in Fig. \ref{fig:3a}. Naturally arises the question if the
differential equation \eqref{exa:elsotetel2D} can be represented
with a mechanism only containing three complexes \ce{2X}, \ce{2Y},
\ce{X + Y}. It can be easily shown that the answer is negative.
\end{Exa}

\begin{Exa}
Let $M=3$ and suppose that $V(x,y,z)=x^2+y^2+z^2$. Then
\eqref{ered:elsotetel} specializes to
\begin{eqnarray}
{x}'=ay^2+bz^2-cxy-exz\nonumber\\
{y}'=cx^2+dz^2-axy-fyz\label{exa:elsotetel3D}\\
{z}'=ex^2+fy^2-bxz-dyz\nonumber
\end{eqnarray}
(with nonnegative coefficients $a,b,c,d,e,f$) which may be
considered as the \ikde\ of the reaction shown in Fig.
\ref{example2}. as again the application of \mma{RightHandSide}
verifies. A typical trajectory is shown in Fig. \ref{fig:3b}.
\begin{figure}[!ht]
\centering
  \includegraphics[width=0.45\textwidth]{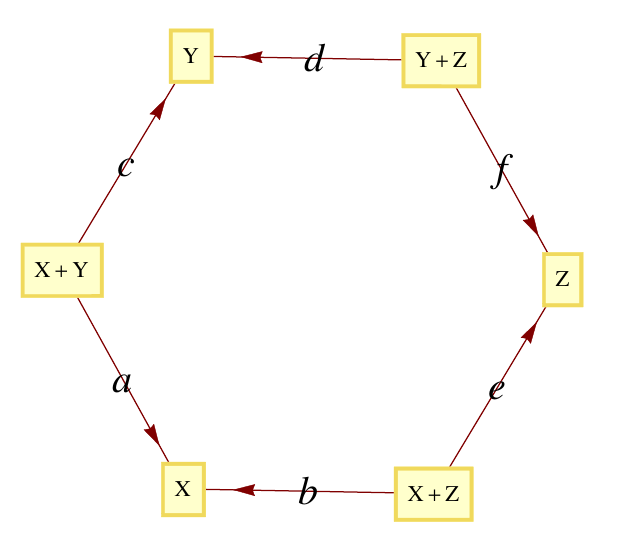}\quad
  \includegraphics[width=0.45\textwidth]{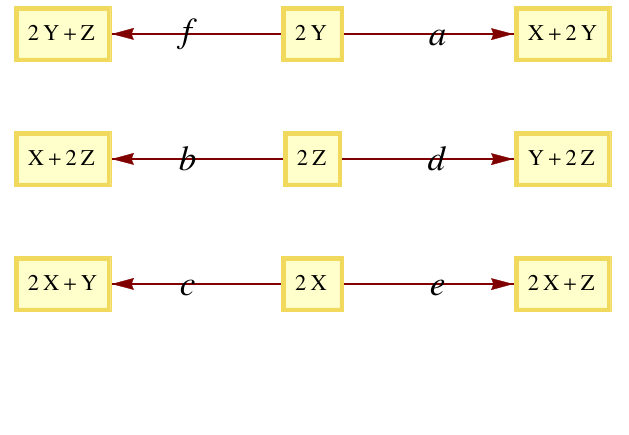}
\caption{3D system with a quadratic first integral} \label{example2}
\end{figure}
\end{Exa}

\begin{Cor}
As the divergence of the system \eqref{exa:elsotetel3D} is
$-ax-bx-cy-dy-ez-fz<0$ in the first orthant and the system has a
first integral, \cite[Theorem 3.3]{tothZAMM} (actually, a version of
K. R. Schneider's theorem) implies that it has no periodic orbit in
the first orthant.
\end{Cor}

The next result shows that (even weighted) sum of squares cannot be
a first integral if mass is conserved.

\begin{Thm}
Let us consider the differential equation system \eqref{system1}
where the functions $F_m$ are of the form \eqref{quadratic1}.
Suppose that the differential equation system is \emph{kinetic} and
\emph{kinetically mass conserving}. The function
$$
V(x_1,x_2,\dots,x_M)=a_1x_1^2+a_2x_2^2+\dots +a_Mx_M^2
$$
(where $a_m \neq 0$ for all $m$) is a first integral for the system
\eqref{system1}, if and only if for all $m:$ $$
F_m(x_1,x_2,\dots,x_M)=0.
$$
\end{Thm}

\emph{Proof}. The function $V$ is a first integral for the system
\eqref{system1} if and only if \eqref{Lie-derivative} holds. Since
the system \eqref{system1} is kinetic and mass conserving with some
positive numbers $\varrho_m, m=1,2,\dots,M$; besides $(k1)$--$(k4)$
the following inequalities also hold:

\noindent
\begin{tabular}{lllll}
$(m1)$&for all $m$& $\displaystyle\sum_{p=1}^M \varrho_pA_{p,m}=0$  (the sum of the coefficients of $x_m^2$)\\

$(m2)$&for all $m,p$& $\displaystyle
\varrho_mB_{m,p}+\varrho_pB_{p,m}
+\sum_{\begin{smallmatrix}q=1\\q\neq m,q\neq
p\\m<p\end{smallmatrix}}^{M}\varrho_qC_{m,p}^q=0$\\&&  (the sum of the coefficients of $x_m x_p$)\\
$(m3)$&for all $m$& $\displaystyle\sum_{p=1}^M \varrho_pD_{p,m}=0$  (the sum of the coefficients of $x_m$)\\
$(m4)$&for all $m$& $\displaystyle\sum_{m=1}^M \varrho_mE_m=0.$
\end{tabular}

Similarly, as in the proof of Theorem \ref{th:elsotetel},  for all
$m: A_{m,m}=0,\;D_{m,m}=0$ and $E_m=0.$ Then,  for all $m\neq p$,
because of $(k1)$ and $(m1): A_{m,p}=0,$ and because of $(k3)$ and
$(m3): D_{m,p}=0$.

If $m\neq p$, then the coefficient of $x_m^2x_p$ in
\eqref{Lie-derivative} is $a_mB_{m,p}+a_pA_{p,m}=0$. Since
$A_{p,m}=0$, it follows that  for all $m,p, m\neq p:$ $B_{p,m}=0$.
Finally,  because of $(k2)$ and $(m2)$, for all $m,p,q$ such that
$q\neq m,q\neq p,m<p: C_{m,p}^q=0$.

The proof of the \emph{if} part is obvious.


\begin{Thm}
Let us consider the following differential equation system
\begin{eqnarray}
x_k'&=&X_k\circ(x_1,\dots,x_K,y_1,\dots,y_L,z), \quad (k=1,\dots,K)\nonumber\\
y_l'&=&Y_l\circ(x_1,\dots,x_K,y_1,\dots,y_L,z), \quad (l=1,\dots,L)\label{system2}\\
z'&=&Z\circ(x_1,\dots,x_K,y_1,\dots,y_L,z)\nonumber
\end{eqnarray}
where the functions $X_k,Y_l,Z$ are quadratic polynomials of all the
variables and $K,L$ are positive integers. Suppose that the
differential equation system is \emph{kinetic} and \emph{kinetically
mass conserving} with some positive numbers $\varrho_k^x,
\varrho_l^y, \varrho^z$. $(k=1,\dots,K, l=1,\dots,L)$. The function
$$
V(x_1,\dots,x_K,y_1,\dots,y_L,z)=a_1x_1^2+\dots
a_Kx_K^2-b_1y_1^2-\dots b_Ly_L^2
$$
(where $a_k > 0$ and $b_l > 0$ for all $k,l$) is a first integral
for system \eqref{system2} if and only if
\begin{eqnarray}
X_k(x_1,\dots,x_K,y_1,\dots,y_L,z)&=&\sum_{l=1}^L b_l A_{k,l} y_l z, \quad (k=1,\dots,K)\nonumber\\
Y_l(x_1,\dots,x_K,y_1,\dots,y_L,z)&=&\sum_{k=1}^K a_k A_{k,l} x_k z, \quad (l=1,\dots,L)\nonumber\\
Z&=&-\left(\sum_{k=1}^K \varrho_k^x X_k+\sum_{l=1}^L \varrho_l^y
Y_l\right)\nonumber
\end{eqnarray}
where $A_{k,l}\geq0$ for all $k,l$. If $K=0$ or $L=0$, then
$X_k=Y_l=Z=0$. \label{th:harmadiktetel}
\end{Thm}

\emph{Proof}. The function $V$ is a first integral for the system
(\ref{system2}) if and only if the Lie-derivative is equal to zero,
that is,
\begin{equation}
\sum_{k=1}^K a_k x_k X_k(x_1,\dots,x_K,y_1,\dots,y_L,z)-\sum_{l=1}^L
b_l y_l Y_l(x_1,\dots,x_K,y_1,\dots,y_L,z)=0 \label{Lie-derivative2}
\end{equation}

It is obvious that the functions $X_k,Y_l,Z$ may not contain
constant terms. Since the system \eqref{system2} is kinetic and mass
conserving, the constants are nonnegative and their weighted sum
with positive weights is equal to zero, and thus each constant is
equal to zero.

Next we show that the terms $x_k^2, y_l^2, z^2$ and $x_k, y_l, z$ in
$X_k,Y_l,Z$ have zero coefficients for all $k,l$. Let us consider at
first only the terms containing quadratic and linear monomials.
\begin{eqnarray}
X_k&=& \left(\sum_{i=1}^K A_{k,i}^x x_i^2 +\sum_{j=1}^L B_{k,j}^x
y_j^2+C_k^x z^2\right)+ \left(\sum_{i=1}^K \overline{A}_{k,i}^{\,x}
x_i +\sum_{j=1}^L \overline{B}_{k,j}^{\,x} y_j
+\overline{C}_k^{\,x} z\right)+\dots\nonumber\\
Y_l&=&\left(\sum_{m=1}^K A_{l,m}^y x_m^2 +\sum_{n=1}^L B_{l,n}^y
y_n^2+C_l^y z^2\right)+ \left(\sum_{m=1}^K \overline{A}_{l,m}^{\,y}
x_m +\sum_{n=1}^L \overline{B}_{l,n}^{\,y} y_n
+\overline{C}_l^{\,y} z\right)+\dots\nonumber\\
Z&=&\left(\sum_{k=1}^K A_{k}^z x_k^2+\sum_{l=1}^L B_{l}^z y_l^2+C^z
z^2\right)+ \left(\sum_{k=1}^K \overline{A}_{k}^{\,z}
x_k+\sum_{l=1}^L \overline{B}_{l}^{\,z} y_l+\overline{C}^{\,z}
z\right)+\dots\label{eleje}
\end{eqnarray}

a) The coefficients of $x_k^3$ and $y_l^3$ in
\eqref{Lie-derivative2} are $a_k A_{k,k}^x$ and $-b_l B_{l,l}^y$,
respectively, and so  for all $k,l: A_{k,k}^x=B_{l,l}^y=0$. Because
of mass conservation, the weighted sums of the coefficients of
$x_k^2$ and $y_l^2$ in \eqref{eleje} are equal to zero, that is,
\begin{eqnarray*}
\sum_{s=1,s\neq k}^K \varrho_s^x A_{s,k}^x +\sum_{t=1}^L \varrho_t^y A_{t,l}^y + \varrho^z A_k^z&=&0 \quad \text{for all } 1\leq k\leq K\\
\sum_{s=1}^K \varrho_s^x B_{s,k}^x +\sum_{t=1,t\neq l}^L \varrho_t^y B_{t,l}^y + \varrho^z B_l^z&=&0 \quad \text{for all } 1\leq l\leq L\\
\end{eqnarray*}
Since the system \eqref{system2} is kinetic, each term in the above
sums is nonnegative, therefore each coefficient of $x_k^2$ and
$y_l^2$ in \eqref{eleje} is equal to zero. The coefficients of $x_k
z^2$ and $y_l z^2$ in \eqref{Lie-derivative2} are $a_k C_{k}^x$ and
$-b_l C_{l}^y$, respectively, and so  for all $k,l:
C_{k}^x=C_{l}^y=0$. Therefore, because of mass conservation, $C^z=0$
as well.

b) The coefficients of $x_k^2$ and $y_l^2$ in
\eqref{Lie-derivative2} are $a_k\overline{A}_{k,k}^{\,x}$ and
$-b_l\overline{B}_{l,l}^{\,y}$, respectively, therefore
$\overline{A}_{k,k}^{\,x}=\overline{B}_{l,l}^{\,y}=0$. It can be
shown very similarly as in part a) that the coefficients of $x_k,
y_l, z$ in $X_k,Y_l,Z$ are equal to zero as well.

c) Next consider the terms of the form $x_s x_t$, $y_s y_t$ and $x_i
y_j$ (where $s<t$) in $X_k,Y_l,Z$. We show that these terms have
zero coefficients as well for all $s,t,i,j$, where $s<t$.

\begin{eqnarray}
X_k&=&\sum_{\begin{smallmatrix}s,t=1\\s< t\end{smallmatrix}}^K
D_{s,t}^{x,k} x_s x_t +\sum_{\begin{smallmatrix}s,t=1\\s<
t\end{smallmatrix}}^L E_{s,t}^{x,k} y_s y_t+
\sum_{\begin{smallmatrix}1\leq s\leq K\\1\leq t\leq
L\end{smallmatrix}}F_{s,t}^{x,k} x_s y_t+\dots
\nonumber\\
Y_l&=&\sum_{\begin{smallmatrix}s,t=1\\s< t\end{smallmatrix}}^K
D_{s,t}^{y,l} x_s x_t +\sum_{\begin{smallmatrix}s,t=1\\s<
t\end{smallmatrix}}^L E_{s,t}^{y,l} y_s y_t+
\sum_{\begin{smallmatrix}1\leq s\leq K\\1\leq t\leq
L\end{smallmatrix}}F_{s,t}^{y,l} x_s y_t+\dots
\nonumber\\
Z&=&\sum_{\begin{smallmatrix}s,t=1\\s< t\end{smallmatrix}}^K
D_{s,t}^{z} x_s x_t +\sum_{\begin{smallmatrix}s,t=1\\s<
t\end{smallmatrix}}^L E_{s,t}^{z} y_s y_t+
\sum_{\begin{smallmatrix}1\leq s\leq K\\1\leq t\leq
L\end{smallmatrix}}F_{s,t}^{z} x_s y_t+\dots
\label{kozepe}
\end{eqnarray}
For all $s<t$ the coefficients of $x_s^2 x_t$, $x_s x_t^2$,
$y_s^2y_t$ and $y_s y_t^2$ in (\ref{Lie-derivative2}) are $a_s
D_{s,t}^{x,s}$, $a_t D_{s,t}^{x,t}$, $-b_sE_{s,t}^{y,s}$ and
$-b_tE_{s,t}^{y,t}$, respectively, and so
$D_{s,t}^{x,s}=D_{s,t}^{x,t}=E_{s,t}^{y,s}=E_{s,t}^{y,t}=0$. Since
the system is kinetic and because of mass conservation, for all
$1\leq k<l\leq K$, $\displaystyle\sum_{\begin{smallmatrix}i=1\\i\neq
k,l\end{smallmatrix}}^K \varrho_i^x D_{k,l}^{x,i}+\sum_{j=1}^L
\varrho_j^y D_{k,l}^{y,j}+\varrho^z D_{k,l}^z=0$ and for all $1\leq
m<n\leq L$, $\displaystyle\sum_{i=1}^K \varrho_i^x
E_{m,n}^{x,i}+\sum_{\begin{smallmatrix}j=1\\j\neq
m,n\end{smallmatrix}}^L \varrho_j^y E_{m,n}^{y,j}+\varrho^z
E_{m,n}^z=0$ where each term in the sums is nonnegative and thus the
coefficients of the terms $x_s x_t$ and $y_s y_t$ ($s<t$) are equal
to zero.

d) The coefficients of $x_k^2 y_l$ and $x_k y_l^2$ ($1\leq k\leq
K,1\leq l\leq L$) in \eqref{Lie-derivative2} are $a_kF_{k,l}^{x,k}$
and $-b_lF_{k,l}^{y,l}$, respectively, and so for all $k,l$:
$F_{k,l}^{x,k}=F_{k,l}^{y,l}=0$. Since the system is kinetic, and
because of mass conservation,
$\displaystyle\sum_{\begin{smallmatrix}i=1\\i\neq
k\end{smallmatrix}}^K \varrho_i^x
F_{k,l}^{x,i}+\sum_{\begin{smallmatrix}j=1\\j\neq
l\end{smallmatrix}}^L \varrho_j^y F_{k,l}^{y,j}+\varrho^z
F_{k,l}^z=0$ where each term in the sum is nonnegative. Therefore,
each coefficient of $x_ky_l$ in \eqref{kozepe} is equal to zero.

Now we may suppose that $X_k,Y_l,Z$ have the following form:

\begin{eqnarray}
X_k&=&\sum_{s=1}^K G_{k,s}^{x} x_s z +\sum_{s=1}^L H_{k,s}^{x} y_s z
\nonumber\\
Y_l&=&\sum_{s=1}^K G_{l,s}^{y} x_s z +\sum_{s=1}^L H_{l,s}^{y} y_s z
\nonumber\\
Z&=&\sum_{s=1}^K G_{s}^{z} x_s z +\sum_{s=1}^L H_{s}^{z} y_s z
\label{vege}
\end{eqnarray}

e) The coefficients of $x_k^2 z$ and $y_l^2 z$ in
\eqref{Lie-derivative2} are $a_k G_{k,k}^x$ and $-b_l H_{l,l}^y$,
respectively, thus $G_{k,k}^x=0$ and $H_{l,l}^y=0$ for all $k,l$
($1\leq k\leq K$, $1\leq l\leq L$). The coefficient of $x_k x_l z$
($k\neq l$) in \eqref{Lie-derivative2} is $a_k G_{k,l}^x+a_l
G_{l,k}^x=0$. Since the system is kinetic, $G_{k,l}^x$ and
$G_{l,k}^x$ are nonnegative and thus these coefficients are equal to
zero for all $k,l$ ($1\leq k\leq K, 1\leq l\leq K, k\neq l$).
Similarly, $H_{j,s}^{y}=0$ for all $1\leq j\leq L,1\leq s\leq L,
j\neq s$.

f) Finally, the coefficients of $x_k y_l z$ in
(\ref{Lie-derivative2}) is $a_k H_{k,l}^x-b_l G_{l,k}^y=0$ where
$H_{k,l}^x$ and $G_{l,k}^y$ are nonnegative for all $k,l$ ($1\leq
k\leq K,1\leq l\leq L$). Whithout the loss of generality, it may be
assumed that $H_{k,l}^x=b_lA_{k,l}$ where $A_{k,l}\geq0$. Thus
$G_{l,k}^y=a_kA_{k,l}$. Using that the system is kinetically mass
conserving, we obtain the formula for $Z$.

g) If $K=0$ or $L=0$, then it can be shown easily that
$X_k=Y_l=Z=0$. If for example $L=0$, then the same proof can be
repeated with $b_l=-c_l$ where $c_l>0$. Then in case f) $a_k
H_{k,l}^x+c_l G_{l,k}^y=0$ where $a_k>0,c_l>0, H_{k,l}^x\geq0,
G_{l,k}^y\geq0$ and thus $H_{k,l}^x=G_{l,k}^y=0$.

\begin{Exa}
If $V(x,y,z)=x^2-y^2$ and $\varrho^x=\varrho^y=\varrho^z=1$, then
the equation system is ($a\geq0$):
\begin{eqnarray}
x'&=&ayz\nonumber\\
y'&=&axz\nonumber\\
z'&=&-axz-ayz \label{exa:harmadiktetel1}
\end{eqnarray}
A possible reaction is the following:
\begin{eqnarray*}
&&\ce{X + Z ->[a] X + Y <-[a] Y + Z}
\end{eqnarray*}
\end{Exa}

\begin{figure}[!ht]
\centering
  \includegraphics[width=0.45\textwidth]{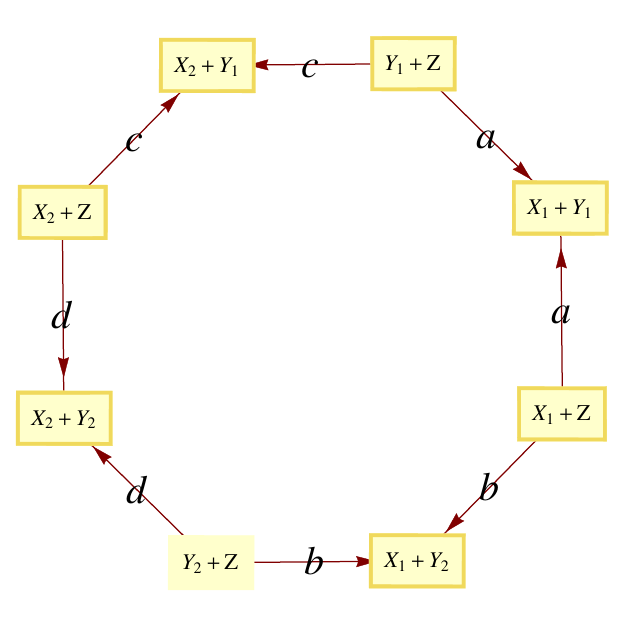}
\caption{Reaction system for \eqref{exa:harmadiktetel2}.}
\label{fig:4}
\end{figure}

\begin{Exa}
If $V(x_1,x_2,y_1,y_2,z)=x_1^2+x_2^2-y_1^2-y_2^2$ and
$\varrho_1^x=\varrho_2^x=\varrho_1^y=\varrho_2^y=\varrho^z=1$, then
the equation system is ($a,b,c,d\geq0$):
\begin{eqnarray}
&&x_1'=ay_1z+by_2z\qquad y_1'=ax_1z+cx_2z\nonumber\\
&&x_2'=cy_1z+dy_2z\qquad y_2'=bx_1z+dx_2z\nonumber\\
&&z'=-x_1'-x_2'-y_1'-y_2' \label{exa:harmadiktetel2}
\end{eqnarray}
A possible reaction is the following:
\begin{eqnarray*}
&&\ce{X_1 + Z ->[1] aY_1 + bY_2 + X_1 + ( 1-a-b)Z}\\
&&\ce{X_2 + Z ->[1] cY_1 + dY_2 + X_2 + ( 1-c-d)Z}\\
&&\ce{Y_1 + Z ->[1] aX_1 + cX_2 + Y_1 + ( 1-a-c)Z}\\
&&\ce{Y_2 + Z ->[1] bX_1 + dX_2 + Y_2 + ( 1-b-d)Z}
\end{eqnarray*}
\end{Exa}
Another possible reaction can be seen in Fig. \ref{fig:4}.

\subsection{The first integral is a binary quadratic form}

Now we investigate first integrals that are quadratic homogeneous
polynomials in two variables, that is, $V(x,y)=ax^2+2bxy+cy^2$.
Obviously, if $V$ is a first integral for a system, then any
nonnegative constant multiples of it is also a first integral for
the same system. Thus, without the loss of generality, it may be
assumed that $a>0$ and $b\neq0$. Consider the following differential
equation system
\begin{alignat}{12}
{x}'&=&A_1x^2&+&B_1xy&+&C_1y^2&+&D_1x&+&E_1y&+&F_1\nonumber\\
{y}'&=&A_2x^2&+&B_2xy&+&C_2y^2&+&D_2x&+&E_2y&+&F_2\label{eq2Dgeneral}
\end{alignat}
and suppose that the differential equation system is \emph{kinetic}.
Then the following statements hold.


\begin{Thm}
The function $V(x,y)=a x^2+2b xy+c y^2$ where $a>0$, $c>0$,
$ac-b^2\neq0$ is a first integral for the system
(\ref{eq2Dgeneral}), if and only if it has the following form
($K\geq0,L\geq0$):
\begin{alignat}{12}
{x}'&=&-bK x^2&+&(-cK+bL)xy&+&cLy^2\nonumber\\
{y}'&=&aK x^2&+&(bK-aL)xy&-&bLy^2\label{eq2D-1}
\end{alignat}
If the system (\ref{eq2D-1}) is kinetically mass conserving then
$x'=y'=0$. \label{Th-eq2D-1}
\end{Thm}


\begin{Thm}
The function $V(x,y)=a x^2+2b xy+c y^2$ where $a>0$, $b>0$, $c>0$,
$ac-b^2=0$ is a first integral for the system (\ref{eq2Dgeneral}),
if and only if it has the following form
($K\geq0,L\geq0,M\geq0,N\geq0$, $S$ is arbitrary):
\begin{alignat}{12}
{x}'&=&-bK x^2&+&cS xy&+&cLy^2&-&bMx&+&cNy\nonumber\\
{y}'&=&aK x^2&-&bS xy&-&bLy^2&+&aMx&-&bNy\label{eq2D-2a}
\end{alignat}
If the system (\ref{eq2D-2a}) is kinetically mass conserving with
some positive numbers $\varrho_1,\varrho_2$ then it has the form
\begin{alignat}{12}
{x}'&=&-\varrho_2 K x^2&+&\varrho_2 S xy&+&\varrho_2 Ly^2&-&\varrho_2 Mx&+&\varrho_2 Ny\nonumber\\
{y}'&=&\varrho_1 K x^2&-&\varrho_1 S xy&-&\varrho_1 Ly^2&+&\varrho_1
Mx&-&\varrho_1 Ny
\end{alignat}
\label{Th-eq2D-2a}
\end{Thm}


\begin{Thm}
The function $V(x,y)=a x^2-2b xy+c y^2$ where $a>0$, $b>0$, $c>0$,
$ac-b^2=0$ is a first integral for the system (\ref{eq2Dgeneral}),
if and only if it has the following form
($K\geq0,L\geq0,M\geq0,N\geq0,R\geq0$, $S$ is arbitrary):
\begin{alignat}{13}
{x}'&=&bK x^2&+&cS xy&+&cLy^2&+&bMx&+&cNy&+&cR\nonumber\\
{y}'&=&aK x^2&+&bS xy&+&bLy^2&+&aMx&+&bNy&+&bR\label{eq2D-2b}
\end{alignat}
If the system (\ref{eq2D-2b}) is kinetically mass conserving then
$x'=y'=0$. \label{Th-eq2D-2b}
\end{Thm}


\begin{Thm}
The function $V(x,y)=a x^2+2b xy-c y^2$ where $a>0$, $c>0$, $b\neq0$
is a first integral for the system (\ref{eq2Dgeneral}), if and only
if it has the following form ($K\geq0,L\geq0,M\geq0$):
\begin{alignat}{13}
{x}'&=&-bK x^2&+&(cK-bL)xy&+&cLy^2&-&bMx&+&cMy\nonumber\\
{y}'&=&aK x^2&+&(bK+aL) xy&+&bLy^2&+&aMx&+&bMy\label{eq2D-3}
\end{alignat}
If the system (\ref{eq2D-3}) is kinetically mass conserving then
$x'=y'=0$. \label{Th-eq2D-3}
\end{Thm}


\begin{Thm}
The function $V(x,y)=a x^2+2b xy$ where $a>0$, $b\neq0$ is a first
integral for the system (\ref{eq2Dgeneral}), if and only if it has
the following form ($K\geq0,M\geq0$, $S$ is arbitrary):
\begin{alignat}{13}
{x}'&=&-bK x^2&-&bS xy&&&-&bMx\nonumber\\
{y}'&=&aK x^2&+&(bK+aS) xy&+&bSy^2&+&aMx&+&bMy\label{eq2D-4}
\end{alignat}
If the system (\ref{eq2D-4}) is kinetically mass conserving then
$x'=y'=0$. \label{Th-eq2D-4}
\end{Thm}


\emph{Proof of Theorem \ref{Th-eq2D-1}}. The function $V$ is a first
integral for the system (\ref{eq2Dgeneral}) if and only if its
Lie-derivative with respect to the system is equal to zero, that is,

\begin{eqnarray}
\frac{1}{2}\left((2ax+2by)x'+(2bx+2cy)y'\right)=0\label{Lie2D}
\end{eqnarray}

Equations $(i)-(ix)$ hold since the coefficients of the following
monomials in \eqref{Lie2D} are equal to zero and $(x)$ holds since
the system kinetic:

\begin{tabular}{lll}
$(i)$       &$x^3$: & $aA_1+bA_2=0$\\
$(ii)$      &$y^3$: & $bC_1+cC_2=0$\\
$(iii)$     &$x^2y$:& $bA_1+aB_1+cA_2+bB_2=0$\\
$(iv)$      &$xy^2$:& $bB_1+aC_1+cB_2+bC_2=0$\\
$(v)$       &$x^2$: & $aD_1+bD_2=0$\\
$(vi)$      &$y^2$: & $bE_1+cE_2=0$\\
$(vii)$     &$xy$:  & $bD_1+aE_1+cD_2+bE_2=0$\\
$(viii)$    &$x$:   & $aF_1+bF_2=0$\\
$(ix)$      &$y$:   & $bF_1+cF_2=0$\\
$(x)$       &       &
$C_1\geq0,E_1\geq0,F_1\geq0,A_2\geq0,D_2\geq0,F_2\geq0$
\end{tabular}
\\\\
Since $A_2,C_1\geq0$, without the loss of generality, it can be
assumed that $A_2=aK$ and $C_1=cL$ where $K,L\geq0$. Thus, because
of $(i)$ and $(ii)$, $A_1=-bK$ and $C_2=-bL$. Substituting these
into $(iii)$ and $(iv)$ gives
\begin{eqnarray*}
aB_1+bB_2=-(ac-b^2)K, \quad bB_1+cB_2=-(ac-b^2)L
\end{eqnarray*}
Since $ac-b^2\neq0$, the unique solution of this equation system is
$B_1=-cK+bL$ and $B_2=bK-aL$.

In $(v)$ and $(vi)$, let $D_2:=aM$ and $E_1:=cN$ where $M,N\geq0$.
Thus, $D_1=-bM$ and $E_2=-bN$. Substituting these into $(vii)$ gives
$(ac-b^2)(M+N)=0$. Since $ac-b^2\neq0$, it follows that $M=N=0$ and
thus $D_1=D_2=E_1=E_2=0$. Finally, the unique solution of the system
$(viii)-(ix)$ is $F_1=F_2=0$.

If, moreover, the system (\ref{eq2D-1}) is kinetically mass
conserving, then there exist positive numbers $\varrho_1$ and
$\varrho_2$ such that the following equalities hold:

\begin{tabular}{lll}
$(xi)$       &$K(-\varrho_1 b+\varrho_2 a)=0$\\
$(xii)$      &$K(-\varrho_1 b+\varrho_2 a)+L(\varrho_1 b-\varrho_2 a)=0$\\
$(xiii)$     &$L(\varrho_1 c-\varrho_2 b)=0$
\end{tabular}

If $K\neq0$ and $L\neq0$ then from $(xi)$ and $(xiii)$ it follows
that $ac-b^2=0$ which is a contradiction. If, for example, $K\neq0$
and $L=0$ then from $(xi)$ and $(xii)$ we obtain the same. Thus,
$K=L=0$.

Theorems \ref{Th-eq2D-2a}-\ref{Th-eq2D-4} can be proved very
similarly.

\begin{figure}[h!]
  \centering
  \subfloat[][$a=2, b=1, c=3, K=1, L=1$]{\label{fig:2a}\includegraphics[width=0.42\textwidth]{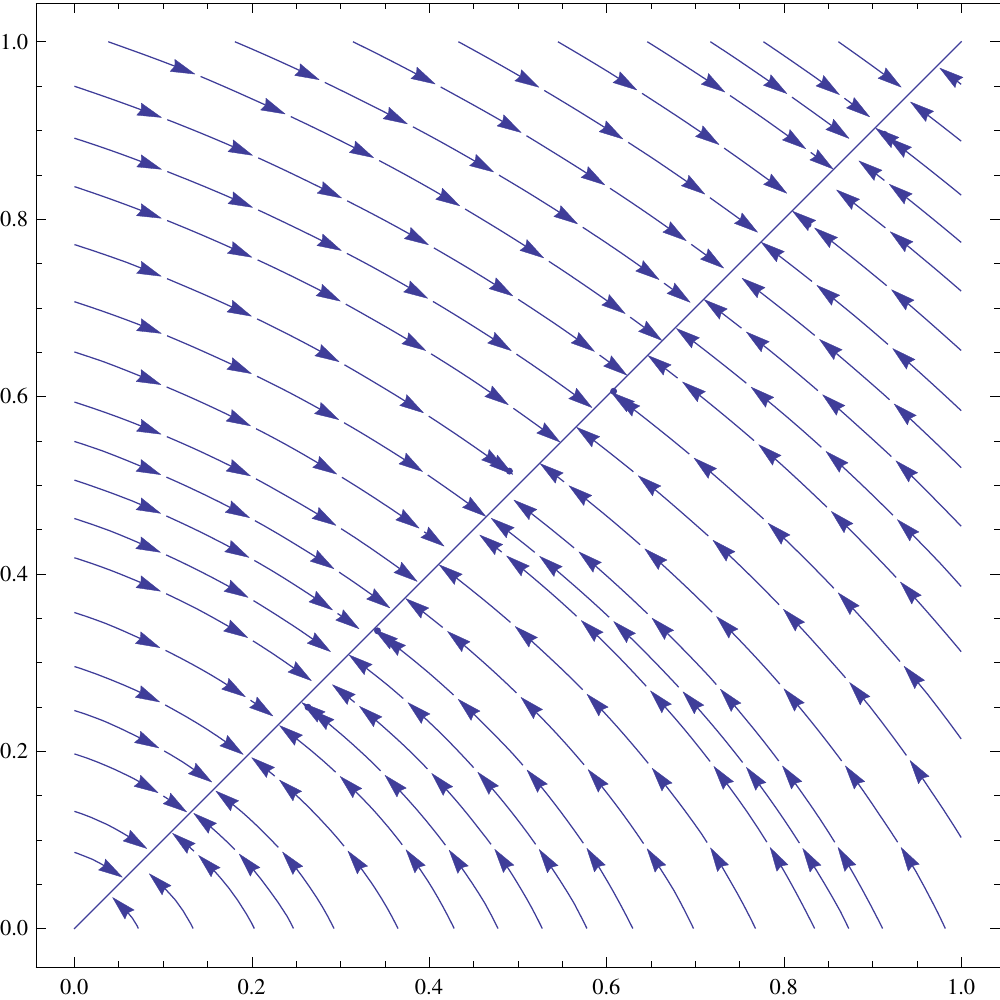}}\qquad
  \subfloat[][$a=1, b=-3, c=2, K=1, L=1$]{\label{fig:2b}\includegraphics[width=0.42\textwidth]{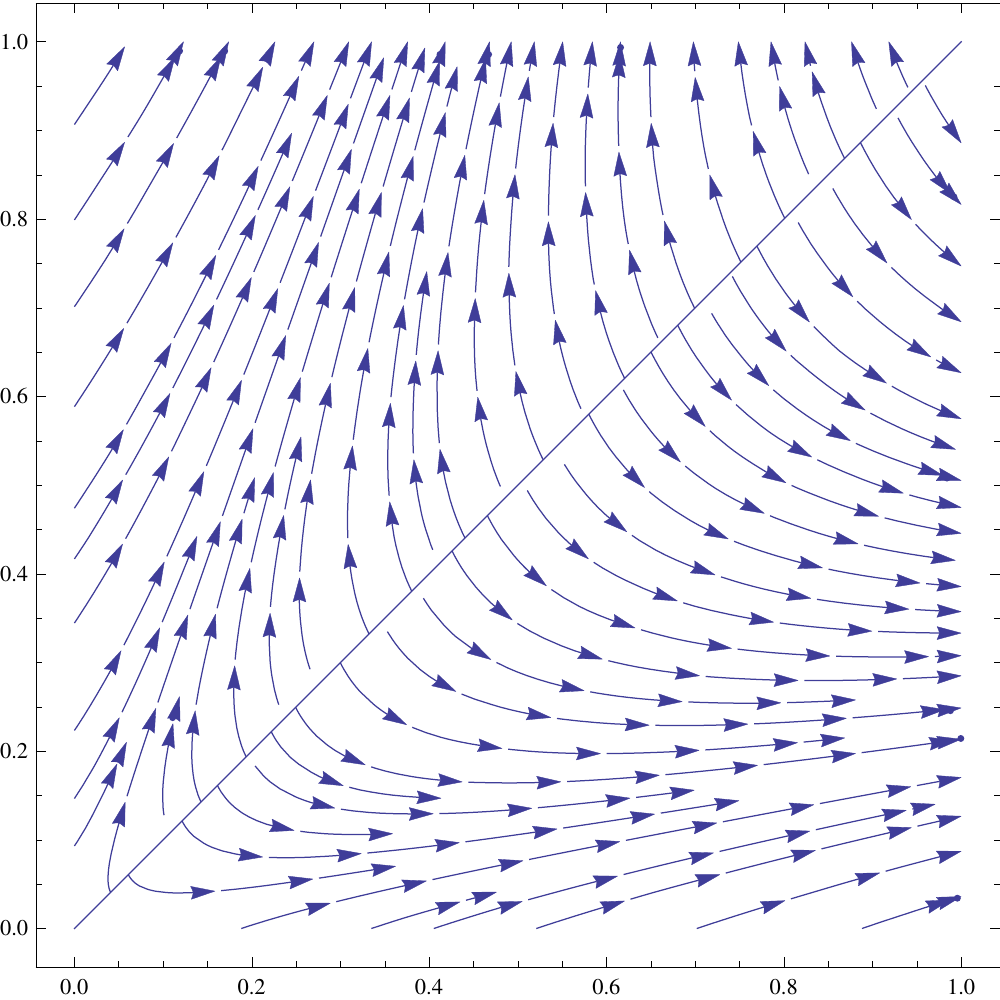}}
  \caption{Trajectories for Example \ref{exa-6} and \ref{exa-7}}\label{fig:2}
\end{figure}

\begin{Rem} The stationary points of system (\ref{eq2D-1}) are (1) $(x,y)=(0,0)$; (2) if $L\neq0$ and $K\neq0$, then $\displaystyle y=\frac{K}{L}x$; (3) if $L=0$, then $x=0$; (4) if $K=0$, then $y=0$. If $ac-b^2>0$, then the trajectories lie on an ellipse while if $ac-b^2<0$, then the trajectories lie on a hyperbola as shown in Fig. \ref{fig:2a} and \ref{fig:2b}. Examples for reactions and trajectories for system (\ref{eq2D-1}) are:

\begin{Exa}
If $a=2, b=1, c=3, K=1, L=1$ then the system is
\begin{alignat}{12}
{x}'&=&-x^2&-&2xy&+&3y^2\nonumber\\
{y}'&=&2x^2&-&xy&-&y^2&\nonumber
\end{alignat}
A possible reaction is the following:
\begin{eqnarray*}
&&\ce{2X + Y <-[2] 2X ->[1] X <-[1] X + Y ->[2] Y <-[1] 2Y ->[3] X +
2Y}
\end{eqnarray*}
\label{exa-6}
\end{Exa}

\begin{Exa}
If $a=2, b=-3, c=2, K=1, L=1$ then the system is
\begin{alignat}{12}
{x}'&=&3x^2&-&5xy&+&2y^2\nonumber\\
{y}'&=&x^2&-&4xy&+&3y^2&\nonumber
\end{alignat}
A possible reaction is the following:
\begin{eqnarray*}
&&\ce{2X + Y <-[1] 2X ->[3] 3X}\\
&&\ce{X + 2Y <-[2] 2Y ->[3] 3Y}\\
&&\ce{X <-[4] X + Y ->[5] Y}
\end{eqnarray*}
\label{exa-7}
\end{Exa}
\end{Rem}

\subsection{Shifted sum of squares}
\begin{Thm}
Let us consider the differential equation system
(\ref{eq2Dgeneral}). Suppose that this system is \emph{kinetic}. The
function $V(x,y)=(x+a)^2+(x+b)^2$ is a first integral for the system
(\ref{eq2Dgeneral}) if and only if it has the following form:
\begin{alignat}{12}
{x}'&=&Ay(y+b)&-&Bx(y+b)&=&Ay^2-Bxy-bBx+bAy\nonumber\\
{y}'&=&Bx(x+a)&-&Ay(x+a)&=&Bx^2-Axy+aBx-aAy
\end{alignat}
where $A\geq0,B\geq0,a\geq0,b\geq0$. If $a<0$, then $B=0$, and if
$b<0$, then $A=0$. From this it follows that there are no periodic
orbits in the first orthant.
\end{Thm}

If $V(x,y)=(x+1)^2+(y+1)^2$, then the equations and the reactions
are

\begin{tabular}{ll}
\makebox(150,150)[l]{\begin{tabular}{l}
${x}'=Ay^2-Bxy-Bx+Ay$\\${y}'=Bx^2-Axy+Bx-Ay$
\end{tabular}}
&
\includegraphics[width=0.31\textwidth]{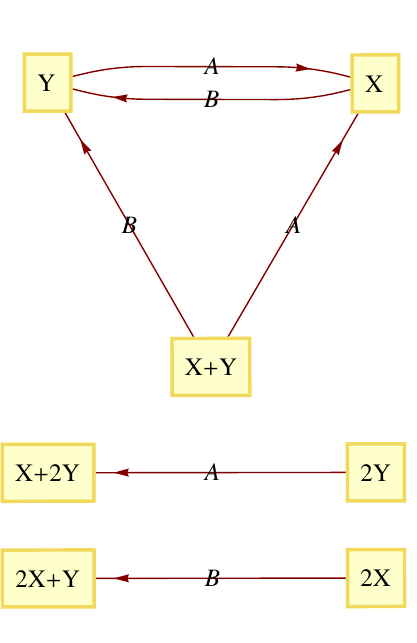}
\end{tabular}

\section{COMPUTER HELP}
Finally let us mention that most of our statements can be obtained
using either the \Mma\ package \mma{ReactionKinetics}
\cite{nagypapptoth,tothnagypapp} or by simple additional programs.
E.g. the following simple code checks if negative cross effect is
present in a polynomial or not.
$$
\mma{CrossEffectQ[polyval\_, vars\_] :=
 Module[\{M=Length[vars]\},}
$$
$$
\mma{And {@}{@} (Map[Not[Negative[\#]] \&,
Flatten[MapThread[ReplaceAll, }
$$
$$
\mma{\{MonomialList[polyval, vars],Thread[vars->\#] \& /{@}}
$$
$$
\mma{(1 - IdentityMatrix[M])\}]]])]}.
$$
And now let us use the newly defined function.
$$
\mma{CrossEffectQ[\{d, c-4yx$^2$+5xy+6z+7w, ax+2y,
-bxy\},\{x,y,z,w\}]}
$$
The answer is as expected depending on the signs of the parameters.
$$
\mma{!Negative[d] \&\& !Negative[c] \&\& !Negative[a] \&\&
!Negative[-b]}
$$
A more easily readable version leads to the same result.
$$
\mma{CrossEffectQ2[polyval\_, vars\_] :=Module[\{M = Length[vars],
L\},}
$$
$$
\mma{L[i\_] := If[Head[polyval[[i]]]
=== Plus,}
$$
$$
\mma{Apply[List, polyval[[i]]], \{polyval[[i]]\}] /.
MapThread[Rule,}
$$
$$
\mma{\{vars, ReplacePart[ConstantArray[1, M], \{i\} -> 0]\}];}
$$
$$
\mma{And @@ Map[\# >= 0 \&, Flatten[Table[L[i], \{i, 1, M\}]]]]}
$$
Using this function for the same example
$$
\mma{CrossEffectQ2[\{d, c-4yx$^2$+5xy+6z+7w, ax+2y,
-bxy\},\{x,y,z,w\}]}
$$
we obtain the following result
$$
\mma{d >= 0 \&\& c >= 0 \&\& a >= 0 \&\& -b >= 0}
$$

\section{DISCUSSION AND OUTLOOK}
\subsection{Other types of first integrals}
We wonder if it is possible to fully characterize those kinetic
differential equations which are of the second degree and have a
general quadratic first integral.

We might find to try other types of first integrals. Let us mention
one simple, still interesting result.
\begin{Sta}
Among the polynomial differential equations of the form
\begin{alignat}{12}\label{second}
{x}'&=&ax^2&+&bxy&+&cy^2&+&dx&+&ey&+&f\nonumber\\
{y}'&=&Ax^2&+&Bxy&+&Cy^2&+&Dx&+&Ey&+&F
\end{alignat}
(defined in the positive quadrant) the only one having
$$V(p,q):=p+q-\ln(p)-\ln(q)$$ as its first integral is
\begin{equation}
{x}'=bxy-bx\quad{y}'=-bxy+by\label{essentially},
\end{equation}
i. e. the Lotka--Volterra equations (allowing possibly time
reversal).
\end{Sta}
Note that the it is not assumed that \eqref{second} is a kinetic
differential equation, and in the result no restriction is made on
the sign of $b.$

This result is very similar to the result leading uniquely to the
Lotka--Volterra model under different circumstances
\cite{morales,hanusse72,hanusse73,tysonlight,pota,tothharsTCA,schumantoth}.

One might try to generalize this result to the multidimensional
case.

Another form of interesting first integrals is a free energy like
function:
$$
V(\cb):=\sum_{m=1}^Mc_m\ln\left(\frac{c_m}{c_m^0}\right),
$$
which turned out to be a useful Lyapunov function for broad classes
of reactions  \cite{hornjackson,volperthudyaev}. Gonzalez-Gascon and Salas
\cite{gonzalezgasconsalas} have systematically found this type of
first integrals (and other types, as well) for three dimensional
Lotka--Volterra systems.

The question arises if these first integrals are kept by some
numerical methods or not.
\subsection{Relations to numerical methods}

Even the simplest kinetic differential equations can only be solved
by numerical methods, therefore the question if such a method is
able to keep important qualitative properties of the models arouse
very early \cite{bertolazzi,farago}. The first positive answers
included numerical methods which keep the property of kinetic
differential equations that starting from a nonnegative initial
concentration they provide solutions which are nonnegative
throughout their total domain of existence
\cite{farago,karatsonkorotov}. Similarly, numerical methods were
constructed to keep linear and quadratic first integrals. The
meaning and existence of linear first integrals have been studied in
detail: they usually represent mass conservation. The existence of a
positive linear first integral together with nonnegativity
\cite{volpert,volperthudyaev} of the solutions implies that the
complete solution of the kinetic differential equation is defined on
the whole real line, which is not necessarily the case for systems
that are not mass conserving. However, quadratic first integrals
were almost neglected.

We have always used global first integrals. Another approach is given by Gonzalez-Gascon and Salas \cite{gonzalezgasconsalas} who started from local first integrals and tried to extend them in cases if it was possible.

\subsection{Acknowledgements}



The present work has partially been supported by the European
Science Foundation Research Networking Programme: Functional
Dynamics in Complex Chemical and Biological Systems,
 by the Hungarian National Scientific Foundation, No. 84060, and
by the COST Action CM901: Detailed Chemical Kinetic Models for
Cleaner Combustion. This work is connected to the scientific program
of the "Development of quality-oriented and harmonized R+D+I
strategy and functional model at BME" project. This project is
supported by the New Sz\'echenyi Plan (Project ID:
T\'AMOP-4.2.1/B-09/1/KMR-2010-0002).
\bibliography{quadratic}
\bibliographystyle{plain}
\end{document}